\documentclass{amsart}\usepackage{amsmath}

\usepackage{amssymb}
\usepackage{graphicx}
\usepackage{amscd}

\newtheorem{theorem}{Theorem}
\theoremstyle{plain}

\numberwithin{equation}{section}

\input{tcilatex}

\begin{document}
\title[Construction of a Family of NAFIL Loops of Odd Order]{\textsf{%
Construction of a Family of NAFIL Loops of Odd Order }$\mathbf{n=2m+1}$}
\author{Raoul E. Cawagas}
\address{Raoul E. Cawagas, SciTech R\&D Center, OVPRD, Polytechnic
University of the Philippines, Manila}
\email{raoulec@pacific.net.ph / raoulec@yahoo.com}
\urladdr{http://www.geocities.com/raoulec2001/raoulweb.htm}
\thanks{\textsl{2000 Mahematics Subject classification}. Primary 20N05;
Secondary 05B15}
\keywords{NAFIL loops, quasigroups, non-associative, Latin squares}

\begin{abstract}
\noindent The existence of NAFIL loops of every odd order $n\geq 5$ is
estab- lished by construction. These are non-associative finite invertible
loops that are simple and power-associative and they form an infinite
family. The first member of this family is the NAFIL loop of order $n=5$
which is known to define a Lie algebra with some possible application in
particle ohysics.
\end{abstract}

\maketitle

\section{Introduction}

In studying any class of finite algebraic structures (like quasigroups,
loops, or groups), the first thing we do is to define the class precisely.
After this, the most important task is to show that the class so defined is
not empty by showing specific examples of its members. However, a class with
only a handful of objects as members is not very interesting. For this
reason, we are more interested in a class with a large number of members.

\emph{Non-associative finite invertible loops (NAFIL)} are loops in which
every element has a unique two-sided inverse and they form an interesting
class that includes the familiar Moufang, Bol, and IP loops$.$ However,
there are other members of this class that have not yet been sufficiently
studied. For instance, several NAFIL loops of small order are known to
define loop algebras that satisfy the Jacobi identity. Some of these loops
are now being studied because of their possible applications in physics \cite%
{ref1}.

\section{On the Existence of NAFIL Loops of Odd Order}

In this paper, we shall prove:

\begin{theorem}
There exists at least one NAFIL loop of every odd order $n\geq 5.$
\end{theorem}

\noindent \textbf{Proof}.

\begin{quote}
To prove this theorem, we shall show how a NAFIL loop $(L_{n},\star )$ of
order $n=2m+1$ can be constructed for any value of $m\geq 2.$ For this, we
need two groups (one of order $m$ and one of order $k=m+1)$ and one
quasigroup of order $k.$

Let $L_{n}=\{1,...,m,m$+$1,...,2m$+$1\}$, where $m\geq 2$, be a set of order 
$n=2m$+$1$ and let $\star $ be a binary operation over $L_{n}$. Next, let $%
L(m)=\{1,...,m\}$ be any group of order $m$ (like the cyclic group $C_{m}$)
and let $L(k)=\{m$+$1,...,2m$+$1\}$ be a group of order $k=m+1$ isomorphic
to the cyclic group $C_{k}$ of order $k$. Hence, $L_{n}=L(m)\cup L(k)$ such
that $L(m)\cap L(k)=\emptyset .$ Moreover, let $\overleftarrow{C}%
_{k}=\{1,...,k\}$ be the counter-cyclic quasigroup \cite{ref2} of order $k$
and let $\overleftarrow{C}_{k}^{T}$ be its transpose.\bigskip

To construct the Cayley table of a system $(L_{n},\star ),$ we proceed as
follows.
\end{quote}

\begin{itemize}
\item First, we form the Latin square blocks $[L(m)],$ $[L(k)]$ (in normal
form) and $[\overleftarrow{C}_{k}]^{T}$ of the systems $L(m),$ $L(k),$ and $%
\overleftarrow{C}_{k}^{T},$ respectively.

\item Second, using the block $[L(k)],$ we form two blocks $[L(k)]^{\prime } 
$ and $[L(k)]^{\prime \prime },$ where
\end{itemize}

\quad (a) $[L(k)]^{\prime }$ is a block of the group $L(k)$ in which row $%
k-1 $ has been omitted.

\quad (b) $[L(k)]^{\prime \prime }$ is a block of the group $L(k)$ in which
column $k$ has been omitted.

\begin{itemize}
\item Third, using the block $[\overleftarrow{C}_{k}]^{T},$ we form another
block $[\overleftarrow{C}_{k}]^{T\ast }$ by replacing each element entry $k$
of $[\overleftarrow{C}_{k}]^{T}$ by elements of the set $L(k)$ as indicated
in Table 2(b).\bigskip
\end{itemize}

\begin{quote}
The block $[L(k)]$ has the general form:\bigskip
\end{quote}

\begin{center}
$\underset{ 
\begin{array}{c}
\\ 
\text{Table 1. General form of the block }[L(k)].%
\end{array}
}{ 
\begin{tabular}{ccccccc}
$m${\small +}$1$ & $m${\small +}$2$ & $m${\small +}$3$ & $\cdots $ & $2m$%
{\small -}$1$ & $2m$ & $2m${\small +}$1$ \\ 
$m${\small +}$2$ & $m${\small +}$3$ & $m${\small +}$4$ & $\cdots $ & $2m$ & $%
2m${\small +}$1$ & $m${\small +}$1$ \\ 
$m${\small +}$3$ & $m${\small +}$4$ & $m${\small +}$5$ & $\cdots $ & $2m$%
{\small +}$1$ & $m${\small +}$1$ & $m${\small +}$2$ \\ 
$\vdots $ & $\vdots $ & $\vdots $ & $\ddots $ & $\vdots $ & $\vdots $ & $%
\vdots $ \\ 
$2m${\small -}$1$ & $2m$ & $2m${\small +}$1$ & $\cdots $ & $2m${\small -}$4$
& $2m${\small -}$3$ & $2m${\small -}$2$ \\ 
$2m$ & $2m${\small +}$1$ & $m${\small +}$1$ & $\cdots $ & $2m${\small -}$3$
& $2m${\small -}$2$ & $2m${\small -}$1$ \\ 
$2m${\small +}$1$ & $m${\small +}$1$ & $m${\small +}$2$ & $\cdots $ & $2m$%
{\small -}$2$ & $2m${\small -}$1$ & $2m$%
\end{tabular}
}$\bigskip
\end{center}

\begin{quote}
Starting with this block $[L(k)]$, we form the blocks $[L(k)]^{\prime }$ (by
deleting row $k-1$ of $[L(k)]$) and $[L(k)]^{\prime \prime }$ (by deleting
column $k$ of $[L(k)]$).\bigskip

Next, we take the block $[\overleftarrow{C}_{k}]^{T}$ of order $k$ which has
the following general form:\bigskip
\end{quote}

\begin{center}
$\underset{ 
\begin{array}{c}
\\ 
\text{Table 2(a). General form of the block }[\overleftarrow{C}_{k}]^{T}%
\text{ of order }k.%
\end{array}
}{ 
\begin{tabular}{ccccccc}
$1$ & \textbf{\emph{k}} & $k$-$1$ & $\cdots $ & $4$ & $3$ & $2$ \\ 
$2$ & $1$ & \textbf{\emph{k}} & $\cdots $ & $5$ & $4$ & $3$ \\ 
$3$ & $2$ & $1$ & $\cdots $ & $6$ & $5$ & $4$ \\ 
$\vdots $ & $\vdots $ & $\vdots $ & $\ddots $ & $\vdots $ & $\vdots $ & $%
\vdots $ \\ 
$k$-$2$ & $k$-$3$ & $k$-$4$ & $\cdots $ & $1$ & \textbf{\emph{k}} & $k$-$1$
\\ 
$k$-$1$ & $k$-$2$ & $k$-$3$ & $\cdots $ & $2$ & $1$ & \textbf{\emph{k}} \\ 
\textbf{\emph{k}} & $k$-$1$ & $k$-$2$ & $\cdots $ & $3$ & $2$ & $1$%
\end{tabular}
}$\bigskip
\end{center}

\begin{quote}
In this block $[\overleftarrow{C}_{k}]^{T}$, replace the entries $k$ by the
entries of the last column of $[L(k)],$ that is, $2m$+$1,m$+$1,...,2m$-$2,2m$%
-$1,2m,$ in this order, from (row 1, column 2), (row 2, column 3) all the
way down to (row $k$-$1$, column $k$), and ending in (row $k,$ column $1)$.
If this is done, we obtain the following block shown in Table 2(b) which we
shall denote by $[\overleftarrow{C}_{k}]^{T\ast }$.\bigskip
\end{quote}

\begin{center}
$\underset{ 
\begin{array}{c}
\\ 
\text{Table 2(b). The block }[\overleftarrow{C}_{k}]^{T\ast }\text{ obtained
from }[\overleftarrow{C}_{k}]^{T}.%
\end{array}
}{ 
\begin{tabular}{ccccccc}
$1$ & \textbf{\emph{2m+1}} & $k$-$1$ & $\cdots $ & $4$ & $3$ & $2$ \\ 
$2$ & $1$ & \textbf{\emph{m+1}} & $\cdots $ & $5$ & $4$ & $3$ \\ 
$3$ & $2$ & $1$ & $\cdots $ & $6$ & $5$ & $4$ \\ 
$\vdots $ & $\vdots $ & $\vdots $ & $\ddots $ & $\vdots $ & $\vdots $ & $%
\vdots $ \\ 
$k$-$2$ & $k$-$3$ & $k$-$4$ & $\cdots $ & $1$ & \textbf{\emph{2m-2}} & $k$-$%
1 $ \\ 
$k$-$1$ & $k$-$2$ & $k$-$3$ & $\cdots $ & $2$ & $1$ & \textbf{\emph{2m-1}}
\\ 
\textbf{\emph{2m}} & $k$-$1$ & $k$-$2$ & $\cdots $ & $3$ & $2$ & $1$%
\end{tabular}
}$
\end{center}

\begin{quote}
Using the blocks $[L(m)],$ $[L(k)]^{\prime },$ $[L(k)]^{\prime \prime },$
and $[\overleftarrow{C}_{k}]^{T\ast }$ thus formed above, we can now
construct the following Cayley table of a NAFIL loop $(L_{n},\star )$ of
order $n=2m+1$ as shown in Table 3(a).\bigskip
\end{quote}

\begin{center}
$\underset{ 
\begin{array}{c}
\\ 
\text{Table 3(a). Cayley table of a NAFIL loop }(L_{n},\star )\text{ of
order }n=2m+1.%
\end{array}
}{ 
\begin{tabular}{|c|ccc|ccc|}
\hline
$\star $ & $1$ & $\cdots $ & $m$ & $m+1$ & $\cdots $ & $2m+1$ \\ \hline
$1$ &  &  &  &  &  &  \\ 
$\vdots $ &  & $[L(m)]$ &  &  & $[L(k)]^{\prime }$ &  \\ 
$m$ &  &  &  &  &  &  \\ \hline
$m+1$ &  &  &  &  &  &  \\ 
$\vdots $ &  & $[L(k)]^{\prime \prime }$ &  &  & $[\overleftarrow{C}%
_{k}]^{T\ast }$ &  \\ 
$2m+1$ &  &  &  &  &  &  \\ \hline
\end{tabular}
}$
\end{center}

\begin{quote}
Since $m$ is finite, the blocks $[L(m)],\;[L(k)]^{\prime },\;[L(k)]^{\prime
\prime },$ and $[\overleftarrow{C}_{k}]^{T\ast }$ can always be constructed
for all values of $m\geq 2$ and $k=m+1$. Thus, the Cayley table of the
system $(L_{n},\star )$ can be constructed for all values of $m\geq 2.$ This
implies that these systems of odd order $n=2m+1$ form an infinite family.

To show that $(L_{n},\star )$ is indeed a NAFIL loop, note that the
resulting block $[L(n)]$ formed by the blocks $[L(m)],\;[L(k)]^{\prime
},\;[L(k)]^{\prime \prime },$ and $[\overleftarrow{C}_{k}]^{T\ast }$ is an $%
n\times n$ Latin square in standard form over the set $L_{n}=%
\{1,...,m,m+1,...,2m+1\},$ that is, the entries in the first row and first
column of $[L(n)]$ are the elements of $L_{n}$ in natural order. If this
block $[L(n)]$ is now converted into the Cayley table shown in Table 3(a),
the element $1$ is seen to be a unique identity element. This means that $%
(L_{n},\star )$ is at least a loop.

Clearly, the group $L(m)$ is a subgroup of $L(n)$. Being a group, every
element of $L(m)$ has a unique inverse. Moreover, every element of the
subset $L(k)$ is of order $2$ (self-inverse) because the square of such an
element is an entry $1$ in the diagonal of the block $[\overleftarrow{C}%
_{k}]^{T\ast }$. Since $L_{n}=L(m)\cup L(k),$ then every element of $L_{n}$
has a unique inverse. Therefore, the system $(L_{n},\star )$ is an
invertible loop.

Finally, the order $m$ of the subgroup $L(m)$ is not a divisor of the order $%
n$ of $(L_{n},\star ).$ By Lagrange's theorem, it follows that $(L_{n},\star
)$ in not a group and hence it is a NAFIL. $\blacksquare $\bigskip
\end{quote}

The smallest NAFIL loop $(L_{5},\star )$ that can be constructed using the
above procedure is of order $n=5$ when $m=2.$ This Cayley table defines a
non-abelian NAFIL loop of order $n=5$ that is simple. Analysis using the
software \emph{FINITAS }\cite{ref3} has shown that it has four subgroups of
order $m=2$ and that it satisfies the cross-inverse property (CIP), the weak
inverse property (WIP), automorphic inverse property (AIP), the flexible law
(FL), power-associative property (PAP), ``A sub m'' loop property (A$_{\text{%
m}}$), and the RIF loop property (RIF).\bigskip

\begin{center}
$\underset{\text{Cayley table of NAFIL loop of order n = 5.}}{\underset{}{ 
\begin{tabular}{|c|cc|ccc|}
\hline
$\star $ & 1 & 2 & 3 & 4 & 5 \\ \hline
1 & 1 & 2 & 3 & 4 & 5 \\ 
2 & 2 & 1 & 5 & 3 & 4 \\ \hline
3 & 3 & 4 & 1 & 5 & 2 \\ 
4 & 4 & 5 & 2 & 1 & 3 \\ 
5 & 5 & 3 & 4 & 2 & 1 \\ \hline
\end{tabular}
}}\medskip $
\end{center}

This loop can be used as the basis of a loop algebra $\mathcal{A}(L_{5})$
whose associated commutator algebra $\mathcal{A}^{-}(L_{5})$ satisfies the
Jacobi identity. Hence $\mathcal{A}^{-}(L_{5})$ is a Lie algebra which has a
subalgebra that is related to the algebra of the Pauli spin matrices in
particle physics.\bigskip

\section{Sample Construction}

\bigskip

We now show how we can construct a loop $(L_{n},\ast )$ of order $n=9$,
where $m=4$ and $k=m+1=5$.\medskip

First, we start with the Latin square blocks $[L(4)],$ $L(5)],$ and $[%
\overleftarrow{C}_{5}]^{T}$ shown below. Second, we form the block $%
[L(5)]^{\prime }$ by deleting row 4 of $[L(5)].$ Third, we form the block $%
[L(5)]^{\prime \prime }$ by deleting column $5$ of $[L(5)].$ And fourth, we
form the block $[\overleftarrow{C}_{5}]^{T\ast }$ by replacing the entries 
\textbf{5} in $[\overleftarrow{C}_{5}]$ by the entries \textbf{9, 5, 6, 7, 8}
in the last row of $[L(5)].$ These blocks are shown below.

\medskip

\begin{center}
\begin{tabular}{ccccccccccccccccccc}
{\small 1} & {\small 2} & {\small 3} & {\small 4} &  &  &  & {\small 5} & 
{\small 6} & {\small 7} & {\small 8} & {\small 9} &  &  & {\small 1} & 
\textbf{5} & {\small 4} & {\small 3} & {\small 2} \\ 
{\small 2} & {\small 3} & {\small 4} & {\small 2} &  &  &  & {\small 6} & 
{\small 7} & {\small 8} & {\small 9} & {\small 5} &  &  & {\small 2} & 
{\small 1} & \textbf{5} & {\small 4} & {\small 3} \\ 
{\small 3} & {\small 4} & {\small 1} & {\small 2} &  &  &  & {\small 7} & 
{\small 8} & {\small 9} & {\small 5} & {\small 6} &  &  & {\small 3} & 
{\small 2} & {\small 1} & \textbf{5} & {\small 4} \\ 
{\small 4} & {\small 1} & {\small 2} & {\small 3} &  &  &  & {\small 8} & 
{\small 9} & {\small 5} & {\small 6} & {\small 7} &  &  & {\small 4} & 
{\small 3} & {\small 2} & {\small 1} & \textbf{5} \\ 
\multicolumn{4}{c}{$\lbrack L(4)]$} &  &  &  & {\small 9} & {\small 5} & 
{\small 6} & {\small 7} & {\small 8} &  &  & \textbf{5} & {\small 4} & 
{\small 3} & {\small 2} & {\small 1} \\ 
&  &  &  &  &  &  & \multicolumn{5}{c}{$[L(5)]$} &  &  & \multicolumn{5}{c}{$%
[\overleftarrow{C}_{5}]^{T}$} \\ 
&  &  &  &  &  &  &  &  &  &  &  &  &  &  &  &  &  &  \\ 
{\small 5} & {\small 6} & {\small 7} & {\small 8} & {\small 9} &  &  & 
{\small 5} & {\small 6} & {\small 7} & {\small 8} &  &  &  & {\small 1} & 
\textbf{9} & {\small 4} & {\small 3} & {\small 2} \\ 
{\small 6} & {\small 7} & {\small 8} & {\small 9} & {\small 5} &  &  & 
{\small 6} & {\small 7} & {\small 8} & {\small 9} &  &  &  & {\small 2} & 
{\small 1} & \textbf{5} & {\small 4} & {\small 3} \\ 
{\small 7} & {\small 8} & {\small 9} & {\small 5} & {\small 6} &  &  & 
{\small 7} & {\small 8} & {\small 9} & {\small 5} &  &  &  & {\small 3} & 
{\small 2} & {\small 1} & \textbf{6} & {\small 4} \\ 
{\small 9} & {\small 5} & {\small 6} & {\small 7} & {\small 8} &  &  & 
{\small 8} & {\small 9} & {\small 5} & {\small 6} &  &  &  & {\small 4} & 
{\small 3} & {\small 2} & {\small 1} & \textbf{7} \\ 
\multicolumn{5}{c}{$\lbrack L(5)]^{\prime }$} &  &  & {\small 9} & {\small 5}
& {\small 6} & {\small 7} &  &  &  & \textbf{8} & {\small 4} & {\small 3} & 
{\small 2} & {\small 1} \\ 
&  &  &  &  &  &  & \multicolumn{4}{c}{$[L(5)]^{\prime \prime }$} &  &  &  & 
\multicolumn{5}{c}{$[\overleftarrow{C}_{5}]^{T\ast }$}%
\end{tabular}

\bigskip
\end{center}

If we now put the blocks $[L(4)],$ $[L(5)]^{\prime },$ $[L(5)]^{\prime
\prime },$ and $[\overleftarrow{C}_{5}]^{T\ast }$ together as indicated in
Table 3(a), we obtain the Cayley table shown in Table 3(b) of a NAFIL loop
of odd order $n=9.$\bigskip

\begin{center}
$\underset{ 
\begin{array}{c}
\\ 
\text{Table 3(b). NAFIL loop}(L_{9},\star )\text{ of order }n=9.%
\end{array}
}{ 
\begin{tabular}{|l|llll|lllll|}
\hline
$\star $ & $1$ & $2$ & $3$ & $4$ & $5$ & $6$ & $7$ & $8$ & $9$ \\ \hline
$1$ & $1$ & $2$ & $3$ & $4$ & $5$ & $6$ & $7$ & $8$ & $9$ \\ 
$2$ & $2$ & $3$ & $4$ & $1$ & $6$ & $7$ & $8$ & $9$ & $5$ \\ 
$3$ & $3$ & $4$ & $1$ & $2$ & $7$ & $8$ & $9$ & $5$ & $6$ \\ 
$4$ & $4$ & $1$ & $2$ & $3$ & $9$ & $5$ & $6$ & $7$ & $8$ \\ \hline
$5$ & $5$ & $6$ & $7$ & $8$ & $1$ & $\mathbf{9}$ & $4$ & $3$ & $2$ \\ 
$6$ & $6$ & $7$ & $8$ & $9$ & $2$ & $1$ & $\mathbf{5}$ & $4$ & $3$ \\ 
$7$ & $7$ & $8$ & $9$ & $5$ & $3$ & $2$ & $1$ & $\mathbf{6}$ & $4$ \\ 
$8$ & $8$ & $9$ & $5$ & $6$ & $4$ & $3$ & $2$ & $1$ & $\mathbf{7}$ \\ 
$9$ & $9$ & $5$ & $6$ & $7$ & $\mathbf{8}$ & $4$ & $3$ & $2$ & $1$ \\ \hline
\end{tabular}
}$\bigskip
\end{center}

This non-abelian NAFIL loop is simple and power-associative. It has one
subgroup of order 4 and six of order 2. \bigskip

\subsection{Some Remarks}

In constructing the NAFIL loop $(L_{n},\star ),$ we made use of the groups $%
L(m)$ and $L(k)$ and the counter-clockwise quasigroup $\overleftarrow{C}%
_{k}^{T}.$ The group $L(m)$ can be any group of order $m$ while the group $%
L(k)=\{m+1,...,2m+1\}$ must be a group of order $k=m+1$ that is isomorphic
to the cyclic group $C_{k}=\{1,...,k]$ of order $k$ whose Latin square block 
$[C_{k}]$ is shown in Table 4.\bigskip

\begin{center}
$\underset{ 
\begin{array}{c}
\\ 
\text{Table 4. Latin square block }[C_{k}]\text{ of the cyclic group }C_{k}.%
\end{array}
}{ 
\begin{tabular}{ccccccc}
1 & 2 & 3 & $\cdots $ & k-2 & k-1 & k \\ 
2 & 3 & 4 & $\cdots $ & k-1 & k & 1 \\ 
3 & 4 & 5 & $\cdots $ & k & 1 & 2 \\ 
$\vdots $ & $\vdots $ & $\vdots $ & $\ddots $ & $\vdots $ & $\vdots $ & $%
\vdots $ \\ 
k-2 & k-1 & k & $\cdots $ & k-5 & k-4 & k-3 \\ 
k-1 & k & 1 & $\cdots $ & k-4 & k-3 & k-2 \\ 
k & 1 & 2 & $\cdots $ & k-3 & k-2 & k-1%
\end{tabular}
}$\bigskip
\end{center}

It is clear that if $(m+i)\in L(k)$ and $i\in C_{k},$ where $i=1,...,k,$
then the following one-to-one correspondence between $L(k)$ and $C_{k}$ is
an isomorphism: $(m+i)\longleftrightarrow i$. Thus, if we simply rename
every element $(m+i)$ of $L(k)]$ by $i,$ we readily obtain the block $%
[C_{k}].$\bigskip

If we permute the rows of $[C_{k}]$ according to the row permutation

\begin{equation*}
_{_{\left( 
\begin{tabular}{ccccccc}
{\small 1} & {\small 2} & {\small 3} & $\cdots $ & {\small k-2} & {\small k-1%
} & {\small k} \\ 
{\small 1} & {\small k} & {\small k-1} & $\cdots $ & {\small 4} & {\small 3}
& {\small 2}%
\end{tabular}
\right) }}
\end{equation*}
then we obtain the block $[\overleftarrow{C}_{k}]^{T}$ shown in Table 2(a)
which is the transpose of the counter-cyclic block $[\overleftarrow{C}_{k}]$.%
\emph{\ }This Latin square block defines a quasigroup of order $k$ denoted
by $\overleftarrow{C}_{k}^{T}$. Since the cyclic group $C_{k}$ exists for
all values of $k\geq 5,$ then $\overleftarrow{C}_{k}^{T}$ also exists for
all values of $k\geq 5.$ Thus, both systems form families.\medskip

It is interesting to note that the quasigroup $\overleftarrow{C}_{k}^{T}$
satisfies the Left Bol (LBol) property.\bigskip

\section{Summary}

\bigskip

In this paper, we proved the existence of at least one NAFIL loop $%
(L_{n},\star )$ of every finite order $n=2m+1,$ where $m\geq 2.$ This was
done by actually constructing $(L_{n},\star )$ using two groups $L(m)$ of
order $m$ and $L(k)$ of order $k=m+1,$ and a special quasigroup $%
\overleftarrow{C}_{k}$ of order $k.$

The group $L(m)$ is any group of order $m$ while $L(k)$ is a group
isomorphic to the cyclic group $C_{k}$ of order $k.$ On the other hand, the
quasigroup $\overleftarrow{C}_{k}$ (which is combinatorially equivalent to $%
C_{k})$ satisfies the Left Bol Property. Using these systems, the NAFIL loop 
$(L_{n},\star )$ is constructed and shown to exist for all values of $%
n=2m+1, $ where $m\geq 2.$ Thus, these loops form an infinite family.

We also indicated that the loop $(L_{5},\star )$ can be used as the basis of
a loop algebra whose associated commutator algebra is a Lie algebra with
potential applications in particle physics.


\begin{thebibliography}{9}
\bibitem{ref1} P. H. Frampton, S. L. Glashow, T. W. Kephart, R. M. Rohm, 
\emph{Non-Associative Loops for Holger Bech Neilsen},
arXiv:hep-th/0111292v1, 30 Nov. 2001

\bibitem{ref2} R. E. Cawagas, \emph{Latin Square Composition of Factorable
Groups and Loops}, Matimyas Matematika, Vol. 21, No. 1, pp. 1-11, (1998)

\bibitem{ref3} R. E. Cawagas, \emph{FINITAS - A Software for the
Construction and Analysis of Finite algebraic Structures}, PUP Journal of
Research and Exposition, Vol. 1, No. 1, 1997.\vspace{0.5in}
\end{thebibliography}
\end{document}